\documentclass[12pt]{article}

\usepackage{amsmath, amssymb}
\usepackage[margin=1in]{geometry}
\usepackage{float}     
\usepackage{placeins}  
\usepackage{tikz}
\usepackage{hyperref}
\usepackage{enumitem}
\usepackage{graphicx}
\usepackage{cite}

\newcommand{\graphscale}{1.4}
\tikzset{
  vertex/.style={circle, fill=black, inner sep=1.5pt},
  every picture/.style={scale=\graphscale, transform shape}
}
\newcommand{\graphbbox}{\path (-1.2,-0.6) rectangle (3.2,3.2);}

\title{Cycles of Length 4 or 8 in Graphs with Diameter 2 and Minimum Degree at Least 3}
\author{Avery Carr\\
Independent Researcher\\
\texttt{avery.carr@ymail.com}
}
\date{Updated: January 29, 2026}

\begin{document}
\maketitle

\begin{abstract}
In this short note it is shown that every graph of diameter $2$ and minimum degree at least $3$
contains a cycle of length $4$ or $8$. This result contributes to the study of the
Erd\H{o}s--Gy\'arf\'as Conjecture \cite{erdos1997} by confirming it for the class of diameter-$2$ graphs.
\end{abstract}

\section*{Notation and Preliminaries}

All graphs considered in this note are finite, simple, and undirected.

Let $G = (V(G), E(G))$ be a graph where $V(G)$ and $E(G)$ denote the set of vertices and edges in G respectively. For a vertex $v \in V(G)$, the \emph{neighborhood} of $v$ is
\[
N(v) = \{ u \in V(G) : uv \in E(G) \},
\]
and the \emph{degree} of $v$ is $d(v) = |N(v)|$. The minimum degree of $G$ is denoted by
\[
\delta(G) = \min\{ d(v) : v \in V(G) \}.
\]

For vertices $u,v \in V(G)$, the \emph{distance} $d(u,v)$ is the length of a shortest path
joining $u$ and $v$ in $G$. The \emph{diameter} of $G$ is
\[
\operatorname{diam}(G) = \max\{ d(u,v) : u,v \in V(G) \}.
\]

A \emph{path} of length $k$ is a sequence of distinct vertices
\[
v_0,v_1,\dots,v_k
\]
such that $v_iv_{i+1} \in E(G)$ for all $0 \le i < k$. A \emph{cycle} of length $k$ is a sequence
\[
v_0,v_1,\dots,v_{k-1},v_0
\]
in which $v_0,\dots,v_{k-1}$ are distinct and $v_iv_{i+1} \in E(G)$ for all indices taken
modulo $k$.

A \emph{$k$-cycle} is a cycle of length $k$. A cycle is called \emph{simple} if it contains no
repeated vertices except for the initial and terminal vertex.

Given a cycle $C$, an edge joining two nonconsecutive vertices of $C$ is called a
\emph{chord} of $C$.

Throughout the paper, vertex labels $v_1,v_2,\dots$ are introduced as needed and are
proved to be distinct. When a sequence is written as
\[
v_1 - v_2 - \cdots - v_k,
\]
it means that $v_iv_{i+1} \in E(G)$ for all $1 \le i < k$.

All set-theoretic notation is used in its standard meaning.

\section*{Introduction}

A well-known open problem of Erd\H{o}s and Gy\'arf\'as asks for unavoidable cycle lengths in graphs with
minimum degree at least three. In particular, they conjectured that every
graph $G$ with $\delta(G)\ge 3$ contains a simple cycle whose length is a power of two.
This is now commonly referred to as the \emph{Erd\H{o}s--Gy\'arf\'as Conjecture}.
Folklore has the conjecture first appearing at a conference in (1995) (later in literature in (1997) \cite{erdos1997}), and is listed in several open-problem compilations
(e.g.\ West \cite{WestOpenProblems} and Erd\H{o}s problems forums such as \cite{ErdosProblems64}).

Despite its simple statement, the conjecture remains open in full generality and has been
verified only for restricted classes of graphs. Early progress includes results for planar
and cubic claw-free graphs \cite{DanielShauger2001, NowbandeganiEtAl2014}.  More recently, Heckman and
Krakovski proved the conjecture for $3$-connected cubic planar graphs \cite{HeckmanKrakovski2013},
and a number of papers have established the conjecture for other hereditary or structured
graph classes such as $P_8$-free graphs \cite{GaoShan2022}. 

Computational searches have also shaped the current understanding. In particular, exhaustive
and heuristic searches (notably by Royle and Markstr\"om) indicate that any counterexample
must be relatively large, and Markstr\"om produced extremal examples illustrating how power-of-two
cycles can be forced to occur only at larger lengths (e.g.\ $16$) in certain cubic graphs
\cite{Markstrom2004}. These computations suggest that minimal counterexamples,
if they exist, are highly constrained.

In this note the Erd\H{o}s--Gy\'arf\'as conjecture is verified for graphs of diameter $2$.
The main result shows that every graph $G$ with $\operatorname{diam}(G)=2$ and $\delta(G)\ge 3$
contains a cycle of length $4$ or $8$ (see Theorem 1.1.). Thus, within the diameter--$2$ regime, the conjecture
holds in its strongest possible form, guaranteeing one of the two smallest nontrivial powers
of two. From the perspective of the broader conjecture, diameter $2$ graphs form a natural
and widely studied class: the global constraint $\operatorname{diam}(G)=2$ forces any two
nonadjacent vertices to have a common neighbor, creating a dense web of short connections.
The proof exploits this interaction between a local degree lower bound and the global
diameter constraint to force short power-of-two cycles.

\newpage

\section*{Main Result}

\textbf{Theorem 1.1.:}
Let $G$ be a graph with diameter $2$ and minimum degree at least $3$.
Then $G$ contains a cycle of length $4$ or $8$.

\section*{Proof}

Assume $G$ has diameter $2$, minimum degree at least $3$, and no $4$-cycle.
For a vertex $v$ in $G$, let $N(v)$ denote its neighborhood; the set of vertices adjacent to $v$.

Let $v_1v_2 \in E(G)$. By the degree condition, $v_1$ has two neighbors other than $v_2$, call them $v_3,v_4$; similarly, $v_2$ has two neighbors $v_5,v_6$. Denote this subgraph of $G$ as $G'$ (see Figure 1).

\begin{figure}[H]
\centering
\begin{tikzpicture}[every node/.style={vertex}]
  \graphbbox
  \node (v1) at (0,0) [label=below:$v_1$] {};
  \node (v2) at (2,0) [label=below:$v_2$] {};
  \node (v3)  at (-0.5,1) [label=above:$v_3$] {};
  \node (v4)  at (0.5,1)  [label=above:$v_4$] {};
  \node (v5)  at (1.5,1)  [label=above:$v_5$] {};
  \node (v6)  at (2.5,1)  [label=above:$v_6$] {};
  \draw (v1)--(v2);
  \draw (v1)--(v3);
  \draw (v1)--(v4);
  \draw (v2)--(v5);
  \draw (v2)--(v6);
\end{tikzpicture}
\caption{$G'$ - Initial edge with neighbors satisfying the degree constraint.}
\label{fig:initial}
\end{figure}
\FloatBarrier

If $v_3 = v_5$ while $v_4 = v_6$, or $v_3 = v_6$ while $v_4 = v_5$, then a cycle of length $4$ forms immediately, namely
\[
v_1-v_3-v_2-v_4-v_1
\quad \text{or} \quad
v_1-v_4-v_2-v_3-v_1.
\]
Both cases contradict the assumption that $G$ contains no $4$-cycle.

Thus, without loss of generality, it suffices to prove the theorem by considering separately
the cases $v_3 = v_5$ and $v_3 \ne v_5$.

\subsection*{Case 1: $v_3 = v_5$}

Assume $v_3 = v_5$. Let $v_4$ be the neighbor of $v_1$ distinct from $v_2,v_3$, and let
$v_6$ be the neighbor of $v_2$ distinct from $v_1,v_3$.  
Set $V(G')=\{v_1,v_2,v_3,v_4,v_6\}$.

\medskip
\noindent\textbf{Claim 1.1.}
$v_4v_6\notin E(G)$.

\emph{Proof.}
If $v_4v_6\in E(G)$, then
\[
v_4 - v_1 - v_2 - v_6 - v_4
\]
is a $4$-cycle, contradicting the assumption that $G$ contains no $4$-cycle.
\hfill$\diamond$

\medskip
Since $v_4$ and $v_6$ are nonadjacent and $\operatorname{diam}(G)=2$, we have
\[
N(v_4)\cap N(v_6)\neq\varnothing.
\]
Choose
\[
v_7\in N(v_4)\cap N(v_6).
\]

\medskip
\noindent\textbf{Claim 1.2.}
$v_7\notin\{v_1,v_2,v_3\}$.

\emph{Proof.}
If $v_7= x$ such that $x \in \{v_1, v_2, v_3\}$, then at least one of the following $4$--cycles would form:
\[
\begin{aligned}
v_1-v_3-v_2-v_6-v_1,\\
v_2-v_3-v_1-v_6-v_2,\\
v_4-v_3-v_2-v_1-v_4.\\
\end{aligned}
\]
\hfill$\diamond$

\medskip
Thus, assume $v_7\notin V(G')$.  
Then the closed walk
\[
v_7-v_4-v_1-v_3-v_2-v_6-v_7
\]
forms a $6$-cycle with a chord $v_1v_2$ (see Figure 2).

\begin{figure}[H]
\centering
\begin{tikzpicture}[every node/.style={vertex}]
  \graphbbox
  \node (v7) at (1,2.4)   [label=above:$v_7$] {};
  \node (v4) at (0,1.3)   [label=left:$v_4$]  {};
  \node (v1) at (0.3,0)   [label=below:$v_1$] {};
  \node (v2) at (1.7,0)   [label=below:$v_2$] {};
  \node (v6) at (2,1.3)   [label=right:$v_6$] {};
  \node (v3) at (1,1) [label=above:$v_3$] {};
  \draw (v7)--(v4)--(v1)--(v2)--(v6)--(v7);
  \draw (v3)--(v1);
  \draw (v3)--(v2);
\end{tikzpicture}
\caption{$6$-Cycle with a $v_1v_2$ Chord}
\label{fig:case1-C6}
\end{figure}

\medskip
Notice $v_3$ and $v_7$ are not adjacent, and if the edge $v_3v_7 \in E(G) $, then a $4$-cycle forms by
\[
v_7-v_3-v_1-v_4-v_7,
\]
a contradiction. Also, the diameter--$2$ condition implies
\[
N(v_3)\cap N(v_7)\neq\varnothing.
\]
Let $v_8\in N(v_3)\cap N(v_7)$.

\medskip
\noindent\textbf{Claim 1.3} $v_8\notin V(G')$.\\
\emph{Proof.}
\smallskip
Indeed, if $v_8\in V(G')$, then $v_8 \in \{v_1, v_2, v_4, v_6\}$ and a $4$--cycle arises in
each case:
\[
\begin{aligned}
v_8=v_1 &\implies v_1-v_2-v_6-v_7-v_1,\\
v_8=v_2 &\implies v_2-v_1-v_4-v_7-v_2,\\
v_8=v_4 &\implies v_4-v_3-v_2-v_1-v_4,\\
v_8=v_6 &\implies v_6-v_3-v_1-v_2-v_6.
\end{aligned}
\]
This contradicts the assumption that $G$ contains no $4$--cycle. Hence
$v_8\notin V(G')$ \\(see Figure 3).
\hfill$\diamond$

\begin{figure}[H]
\centering
\begin{tikzpicture}[every node/.style={vertex}]
  \graphbbox
  \node (v7) at (1,2.4)   [label=above:$v_7$] {};
  \node (v4) at (0,1.3)   [label=left:$v_4$]  {};
  \node (v1) at (0.3,0)   [label=below:$v_1$] {};
  \node (v2) at (1.7,0)   [label=below:$v_2$] {};
  \node (v6) at (2,1.3)   [label=right:$v_6$] {};
  \node (v3) at (1,1) [label=above:$v_3$] {};
  \node (v8) at (2.5,2) [label=above:$v_8$] {};
  \draw (v7)--(v4)--(v1)--(v2)--(v6)--(v7);
  \draw (v7)--(v8);
  \draw (v3)--(v8);
  \draw (v3)--(v1);
  \draw (v3)--(v2);
\end{tikzpicture}
\caption{$v_8 \in N(v_7) \cap N(v_3)$}
\label{fig:case1-v8}
\end{figure}

\medskip
Since $\delta(G)\ge 3$, the vertex $v_8$ has a neighbor
$v_9\notin\{v_3,v_7\}$.

\medskip
\noindent\textbf{Claim 1.4.} $v_9\notin V(G')$.\\
\emph{Proof.} Otherwise $v_9\in\{v_1,v_2,v_4,v_6\}$, and each possibility yields a
$4$--cycle:
\[
\begin{aligned}
v_9=v_1 &\implies v_1-v_4-v_7-v_8-v_1,\\
v_9=v_2 &\implies v_2-v_1-v_3-v_8-v_2,\\
v_9=v_4 &\implies v_4-v_8-v_3-v_1-v_4,\\
v_9=v_6 &\implies v_6-v_2-v_3-v_8-v_6.
\end{aligned}
\]
\hfill$\diamond$

Thus $v_9\notin V(G')$ (see Figure 4).

\begin{figure}[H]
\centering
\begin{tikzpicture}[every node/.style={vertex}]
  \graphbbox
  \node (v7) at (1,2.4)   [label=above:$v_7$] {};
  \node (v4) at (0,1.3)   [label=left:$v_4$]  {};
  \node (v1) at (0.3,0)   [label=below:$v_1$] {};
  \node (v2) at (1.7,0)   [label=below:$v_2$] {};
  \node (v6) at (2,1.3)   [label=right:$v_6$] {};
  \node (v3) at (1,1) [label=above:$v_3$] {};
  \node (v8) at (2.5,2) [label=above:$v_8$] {};
  \node (v9) at (3.5,2) [label=above:$v_9$] {};
  \draw (v7)--(v4)--(v1)--(v2)--(v6)--(v7);
  \draw (v7)--(v8);
  \draw (v8)--(v9);
  \draw (v3)--(v8);
  \draw (v3)--(v1);
  \draw (v3)--(v2);
\end{tikzpicture}
\caption{$v_8v_9 \in E(G)$}
\label{fig:case1-v9}
\end{figure}

\medskip
\noindent\textbf{Claim 1.5.}
If $v_9$ is adjacent to $x \in \{v_1, v_2, v_4, v_6\}$, then there exists a $4$--cycle.\\ 
\emph{Proof.}
\smallskip
Assume $v_9$ is adjacent to $x \in \{v_1, v_2, v_4, v_6\}$. Thus, a $4$--cycle is formed in each of the following cases:
\[
\begin{aligned}
x =v_1 &\implies v_1-v_3-v_8-v_9-v_1,\\
x =v_2 &\implies v_2-v_3-v_8-v_9-v_2,\\
x =v_4 &\implies v_4-v_7-v_8-v_9-v_4,\\
x = v_6 &\implies v_6-v_7-v_8-v_9-v_6.
\end{aligned}
\]
\hfill$\diamond$

Also, since $v_9$ is nonadjacent with both $v_4$ and $v_2$, by $\operatorname{diam}(G)=2$ both $N(v_9)\cap N(v_4)$ and $N(v_9)\cap N(v_2)$ are non-empty.  

\medskip
\noindent\textbf{Claim 1.6.} There is a vertex $v_{10} \in (N(v_9) \cap N(v_4)) \cup (N(v_9) \cap N(v_2))$ such that \\
$v_{10} \notin (\{v_7, v_8\} \cup V(G'))$.\\
\emph{Proof.}
If $v_8 \in (N(v_9) \cap N(v_4)) \cup (N(v_9) \cap N(v_2))$ then $v_8$ is adjacent to one of $v_2$ or $v_4$ and, by Claim 1.4, a $4$--cycle is present.  Thus, $v_{10} \neq v_8$.  Now, suppose $v_{10} \in (N(v_9) \cap N(v_4)) \cup (N(v_9) \cap N(v_2))$ such that $(N(v_9) \cap N(v_4)) \cup (N(v_9) \cap N(v_2)) \subseteq (\{v_7\} \cup V(G'))$. Then, by the diameter 2 condition and Claim 1.5, $v_9$ is adjacent to both $v_3$ and $v_7$ (else $v_9$ would be adjacent to another pair of vertices in $\{v_7\} \cup V(G')$ violating Claim 1.5). However, this forms a $4$--cycle by $v_3-v_8-v_7-v_9-v_3$, providing a contradiction on the claim of no $4$--cycles. 
\hfill$\diamond$

Thus, by Claim 1.6, there is a distinct vertex $v_{10}$ that is in at least one of  $N(v_9) \cap N(v_2)$ or $N(v_9) \cap N(v_4)$ such that $v_{10} \notin (\{v_7, v_8\} \cup V(G'))$.  In either case, an $8$--cycle is formed (see Figure 5 and 6) by: 
\[
\begin{aligned}
v_{10} \in N(v_9) \cap N(v_2) &\implies v_{10}-v_2-v_3-v_1-v_4-v_7-v_8-v_9-v_{10},\\
v_{10} \in N(v_9) \cap N(v_4) &\implies v_{10}-v_4-v_1-v_2-v_6-v_7-v_8-v_9-v_{10}.\\
\end{aligned}
\]

\begin{figure}[H]
\centering
\begin{tikzpicture}[every node/.style={vertex}]
  \graphbbox
  \node (v7) at (1,2.4)   [label=above:$v_7$] {};
  \node (v4) at (0,1.3)   [label=left:$v_4$]  {};
  \node (v1) at (0.3,0)   [label=below:$v_1$] {};
  \node (v2) at (1.7,0)   [label=below:$v_2$] {};
  \node (v6) at (2,1.3)   [label=right:$v_6$] {};
  \node (v3) at (1,1) [label=above:$v_3$] {};
  \node (v8) at (2.5,2) [label=above:$v_8$] {};
  \node (v9) at (3.5,2) [label=above:$v_9$] {};
  \node (v10) at (3,1) [label=below:$v_{10}$] {};
  \draw (v7)--(v4)--(v1)--(v2)--(v6)--(v7);
  \draw (v9)--(v10);
  \draw (v2)--(v10);
  \draw (v7)--(v8);
  \draw (v8)--(v9);
  \draw (v3)--(v8);
  \draw (v3)--(v1);
  \draw (v3)--(v2);
\end{tikzpicture}
\caption{$v_{10} \in N(v_2) \cap N(v_9)$ forming an $8$-Cycle}
\label{fig:case1-v10a}
\end{figure}

\begin{figure}[H]
\centering
\begin{tikzpicture}[every node/.style={vertex}]
  \graphbbox
  \node (v7) at (1,2.4)   [label=above:$v_7$] {};
  \node (v4) at (0,1.3)   [label=left:$v_4$]  {};
  \node (v1) at (0.3,0)   [label=below:$v_1$] {};
  \node (v2) at (1.7,0)   [label=below:$v_2$] {};
  \node (v6) at (2,1.3)   [label=right:$v_6$] {};
  \node (v3) at (1,1) [label=above:$v_3$] {};
  \node (v8) at (2.5,2) [label=below:$v_8$] {};
  \node (v9) at (3.5,2) [label=below:$v_9$] {};
  \node (v10) at (.75,3.25) [label=above:$v_{10}$] {};
  \draw (v7)--(v4)--(v1)--(v2)--(v6)--(v7);
  \draw (v4)--(v10);
  \draw (v9)--(v10);
  \draw (v7)--(v8);
  \draw (v8)--(v9);
  \draw (v3)--(v8);
  \draw (v3)--(v1);
  \draw (v3)--(v2);
\end{tikzpicture}
\caption{$v_{10} \in N(v_4) \cap N(v_9)$ forming an $8$-Cycle}
\label{fig:case1-v10b}
\end{figure}

\subsection*{Case 2: $v_3\neq v_5$}

Assume $v_3\neq v_5$.

\medskip
\noindent\textbf{Claim 2.1.}
If $xy\in E(G)$ for some $x\in\{v_3,v_4\}$ and $y\in\{v_5,v_6\}$, then $G$ contains a $4$-cycle.

\emph{Proof.}
Since $v_1x\in E(G)$ and $v_2y\in E(G)$, the cycle $v_1-x-y-v_2-v_1$ has length $4$.
\hfill$\diamond$

\medskip
Hence we may assume that no edge joins the sets $\{v_3,v_4\}$ and $\{v_5,v_6\}$.

\medskip
\noindent\textbf{Claim 2.2.}
There exist $a\in\{v_3,v_4\}$ and $b\in\{v_5,v_6\}$ such that
\[
N(a)\cap N(b)\not\subseteq \{v_1,v_2\}.
\]

\emph{Proof.}
Because $\operatorname{diam}(G)=2$ and $ab\notin E(G)$ (by Claim 2.1),
every such pair $(a,b)$ has a common neighbor, so $N(a)\cap N(b)\neq\varnothing$.

Suppose for contradiction that for every $a\in\{v_3,v_4\}$ and $b\in\{v_5,v_6\}$ we have
\[
N(a)\cap N(b)\subseteq\{v_1,v_2\}.
\]
Since $v_1$ is adjacent to both $v_3$ and $v_4$, the vertex $v_1$ is a common neighbor of
$(a,b)$ exactly when $v_1b\in E(G)$.
Similarly, since $v_2$ is adjacent to both $v_5$ and $v_6$, the vertex $v_2$ is a common neighbor of
$(a,b)$ exactly when $v_2a\in E(G)$.
Thus, for each pair $(a,b)$, at least one of the edges $v_1b$ or $v_2a$ must exist.

If $v_1$ is adjacent to both $v_5$ and $v_6$, then $v_1-v_5-v_2-v_6-v_1$ is a $4$-cycle.
If $v_2$ is adjacent to both $v_3$ and $v_4$, then $v_1-v_3-v_2-v_4-v_1$ is a $4$-cycle.
Hence $v_1$ is adjacent to at most one of $\{v_5,v_6\}$ and $v_2$ is adjacent to at most one of $\{v_3,v_4\}$.
But then there is some pair $(a,b)$ with $v_1b\notin E(G)$ and $v_2a\notin E(G)$, contradicting the
requirement that each pair $(a,b)$ has a common neighbor in $\{v_1,v_2\}$.
Therefore the claim holds.
\hfill$\diamond$

\medskip
By Claim~2.2, choose $a\in\{v_3,v_4\}$ and $b\in\{v_5,v_6\}$ and a vertex
\[
v_7\in N(a)\cap N(b)\setminus\{v_1,v_2\}.
\]
If $v_7\in\{v_3,v_4,v_5,v_6\}$, then a $4$-cycle occurs:
\[
\begin{aligned}
v_7 = v_3 &\Rightarrow v_1 - v_3 - b - v_2 - v_1,\\
v_7 = v_4 &\Rightarrow v_1 - v_4 - b - v_2 - v_1,\\
v_7 = v_5 &\Rightarrow v_1 - a - v_5 - v_2 - v_1,\\
v_7 = v_6 &\Rightarrow v_1 - a - v_6 - v_2 - v_1.
\end{aligned}
\]

Thus,
\[
v_7\notin V(G')=\{v_1,v_2,v_3,v_4,v_5,v_6\}.
\]
Without loss of generality, let $a=v_3$ and $b=v_5$; hence
\[
v_7\in N(v_3)\cap N(v_5),\quad v_7\notin V(G') 
\]
(see Figure 7).

\begin{figure}[H]
\centering
\begin{tikzpicture}[every node/.style={vertex}]
  \graphbbox
  \node (v1) at (0,0) [label=below:$v_1$] {};
  \node (v2) at (2,0) [label=below:$v_2$] {};
  \node (v3) at (-0.5,1) [label=above:$v_3$] {};
  \node (v4) at (0.5,1)  [label=above:$v_4$] {};
  \node (v5) at (1.5,1)  [label=above:$v_5$] {};
  \node (v6) at (2.5,1)  [label=above:$v_6$] {};
  \node (v7) at (0.5,2)  [label=above:$v_7$] {};
  \draw (v1)--(v2);
  \draw (v1)--(v3);
  \draw (v1)--(v4);
  \draw (v2)--(v5);
  \draw (v2)--(v6);
  \draw (v3)--(v7);
  \draw (v5)--(v7);
\end{tikzpicture}
\caption{A new common neighbor $v_7\in N(v_3)\cap N(v_5)$ with $v_7\notin V(G')$.}
\label{fig:case2-v7}
\end{figure}
\FloatBarrier

\medskip
\noindent\textbf{Claim 2.3.} $v_7v_4\notin E(G)$.

\emph{Proof.}
If $v_7v_4\in E(G)$, then $v_4-v_7-v_3-v_1-v_4$ is a $4$-cycle.
\hfill$\diamond$

\medskip
Consider the nonadjacent pair $(v_4,v_6)$.
If $v_4v_6\in E(G)$, then
Claim 2.1 is violated and there is a $4$-cycle.
Hence $v_4v_6\notin E(G)$ and
\[
N(v_4)\cap N(v_6)\neq\varnothing.
\]
Choose
\[
v_8\in N(v_4)\cap N(v_6).
\]

\medskip
\noindent\textbf{Claim 2.4.}
$v_8\notin\{v_3,v_5,v_7\}$.

\emph{Proof.}
If $v_8=v_3$ or $v_8=v_5$, then Claim~2.1 is violated.
If $v_8=v_7$, then Claim~2.3 is violated.
\hfill$\diamond$

\medskip
\noindent\textbf{Subcase 2A: $v_8=v_1$.}
Then $v_1v_6\in E(G)$. Since $\operatorname{diam}(G)=2$ and $v_6$ and $v_7$ are nonadjacent, we have $$N(v_6) \cap N(v_7) \neq \varnothing.$$

Choose $x \in N(v_6) \cap N(v_7)$. If $x \in V(G') \setminus \{v_6\}$ then a $4$--cycle appears by at least one of the following: 

\[
\begin{aligned}
x = v_1 &\Rightarrow v_1 - v_7 - v_5- v_2- v_1,\\
x = v_2 &\Rightarrow v_2 - v_1 - v_3- v_7- v_2,\\
x = v_3 &\Rightarrow v_3 - v_1 - v_2- v_6- v_3,\\
x = v_4 &\Rightarrow v_4 - v_1 - v_2- v_6- v_4,\\
x = v_5 &\Rightarrow v_5 - v_2 - v_1- v_6- v_5.\\
\end{aligned}
\]

Thus, $x \notin V(G') \setminus \{v_6\}$ (see Figure 8).

\begin{figure}[H]
\centering
\begin{tikzpicture}[every node/.style={vertex}]
  \graphbbox
  \node (v1) at (0,0) [label=below:$v_1$] {};
  \node (v2) at (2,0) [label=below:$v_2$] {};
  \node (v3) at (-0.5,1) [label=above:$v_3$] {};
  \node (v4) at (0.5,1)  [label=above:$v_4$] {};
  \node (v5) at (1.5,1)  [label=above:$v_5$] {};
  \node (v6) at (2.6,1)  [label=above:$v_6$] {};
  \node (v7) at (0.5,2)  [label=above:$v_7$] {};
  \node (x) at (2,2)  [label=above:$x$] {};
  \draw (v1)--(v2);
  \draw (v1)--(v3);
  \draw (v1)--(v4);
  \draw (v2)--(v5);
  \draw (v2)--(v6);
  \draw (v3)--(v7);
  \draw (v5)--(v7);
  \draw (v1)--(v6);
  \draw (x)--(v6);
  \draw (x)--(v7);
\end{tikzpicture}
\caption{Common neighbor $x\in N(v_6)\cap N(v_7)$ with $x\notin V(G') \setminus \{v_6\}$.}
\label{fig:case2-v7}
\end{figure}

Since $v_4$ and $v_5$ are nonadjacent (by Claim 2.1),$$N(v_4) \cap N(v_5) \neq \varnothing,$$ by the $\operatorname{diam}(G)=2$ condition. \\Suppose $y \in N(v_4) \cap N(v_5)$. By Claim 2.1, $y \notin  \{v_3, v_6\}$.  \\ Suppose $y \in \{v_1, v_2, v_7, x\}$. \\

Thus, a $4$--cycle occurs by the following: 

\[
\begin{aligned}
y = v_1 &\Rightarrow v_1 - v_6 - v_2- v_5- v_1,\\
y = v_2 &\Rightarrow v_2 - v_6 - v_1- v_4- v_2,\\
y = v_7 &\Rightarrow v_7 - v_3 - v_1- v_4- v_7,\\
y = x &\Rightarrow x - v_5 - v_2- v_6- x.\\
\end{aligned}
\]

Therefore, $y \in N(v_4) \cap N(v_5)$ such that $y \notin \{v_1, v_2, v_7, x\}$ (see Figure 9) .

\begin{figure}[H]
\centering
\begin{tikzpicture}[every node/.style={vertex}]
  \graphbbox
  \node (v1) at (0,0) [label=below:$v_1$] {};
  \node (v2) at (2,0) [label=below:$v_2$] {};
  \node (v3) at (-0.5,1) [label=above:$v_3$] {};
  \node (v4) at (0.5,1)  [label=above:$v_4$] {};
  \node (v5) at (1.5,1)  [label=above:$v_5$] {};
  \node (v6) at (2.6,1)  [label=above:$v_6$] {};
  \node (v7) at (0.5,2)  [label=above:$v_7$] {};
  \node (x) at (2,2)  [label=above:$x$] {};
  \node (y) at (0.9,1.2)  [label=below:$y$] {};
  \draw (v1)--(v2);
  \draw (v1)--(v3);
  \draw (v1)--(v4);
  \draw (v2)--(v5);
  \draw (v2)--(v6);
  \draw (v3)--(v7);
  \draw (v5)--(v7);
  \draw (v1)--(v6);
  \draw (x)--(v6);
  \draw (x)--(v7);
  \draw (y)--(v5);
  \draw (y)--(v4);
\end{tikzpicture}
\caption{Common neighbor $y\in N(v_4)\cap N(v_5)$ with $y\notin \{v_1, v_2, v_5\}$.}
\label{fig:case2-v7}
\end{figure}

But
\[
x-v_7-v_5-y-v_4-v_1-v_2-v_6-x
\]
forms a cycle of length $8$.\hfill$\diamond$

\medskip
\noindent\textbf{Subcase 2B: $v_8=v_2$.}
This case is handled analogously to $v_8 = v_1$; using diameter and degree constraints, start with $v_2v_4 \in E(G)$, force a $x\in N(v_4) \cap N(v_7)$ and $y \in N(v_3) \cap N(v_6)$, obtaining an $8$--cycle.

\noindent\textbf{Subcase 2C: $v_8\notin\{v_1,v_2\}$.}
Then the edges
\[
v_7v_3,\ v_3v_1,\ v_1v_4,\ v_4v_8,\ v_8v_6,\ v_6v_2,\ v_2v_5,\ v_5v_7
\]
form the simple cycle
\[
v_7-v_3-v_1-v_4-v_8-v_6-v_2-v_5-v_7,
\]
which has length $8$ (see Figure 10).
\newpage
\medskip
Thus, in all cases, $G$ contains a cycle of length $4$ or $8$.

\begin{figure}[H]
\centering
\begin{tikzpicture}[every node/.style={vertex}]
  \graphbbox
  \node (v1) at (0,0) [label=below:$v_1$] {};
  \node (v2) at (2,0) [label=below:$v_2$] {};
  \node (v3)  at (-0.5,1) [label=above:$v_3$] {};
  \node (v4)  at (0.5,1)  [label=above:$v_4$] {};
  \node (v5)  at (1.5,1)  [label=above:$v_5$] {};
  \node (v6)  at (2.5,1)  [label=above:$v_6$] {};
  \node (v7)  at (0.5,2)  [label=above:$v_7$] {};
  \node (v8)  at (1.5,2)  [label=above:$v_8$] {};
  \draw (v1)--(v2);
  \draw (v1)--(v3);
  \draw (v1)--(v4);
  \draw (v2)--(v5);
  \draw (v2)--(v6);
  \draw (v3)--(v7);
  \draw (v5)--(v7);
  \draw (v4)--(v8);
  \draw (v6)--(v8);
\end{tikzpicture}
\caption{Case 2: if $v_3\neq v_5$, an 8-cycle appears.}
\label{fig:case2}
\end{figure}

\hfill$\square$

\section*{Acknowledgments}
The author thanks Dr.\ Michael Albert, Editor-in-Chief of the \textit{Australasian Journal of Combinatorics}, and one anonymous external expert for their careful reading of an earlier draft, advice, and for comments on the originality and merit of this work. Also, thank you to Dr. Tao Wang of Henan University for reading the first draft and presenting a counterexample to case 1 in personal communication that was accounted for in the author's original notes but not in the first draft.

\end{document}